\documentclass{rendiconti}

\RendicontiPagina{1}{Laurent V\'eron}{Boundary value problems relative to heat equation}
\def\DOI{xxxxx} 

\usepackage{amsmath, upgreek}
\usepackage{amssymb}
\usepackage{amscd}
\usepackage{xspace}
\usepackage{verbatim}

\newcommand{\be}{\begin{equation}}
\newcommand{\bel}[1]{\begin{equation}\label{#1}}
\newcommand{\ee}{\end{equation}}
\newcommand{\barr}{\begin{eqnarray}}
\newcommand{\earr}{\end{eqnarray}}
\newcommand{\bars}{\begin{eqnarray*}}
\newcommand{\ears}{\end{eqnarray*}}

\newtheorem{subn}{\name}

\newcommand{\bsn}[1]{\def\name{#1}\begin{subn}}
\newcommand{\esn}{\end{subn}}
\newtheorem{sub}{\name}[section]
\newcommand{\bs}{\begin{sub}}
\newcommand{\es}{\end{sub}}
\newcommand{\bsl}[1]{\begin{sub}\label{#1}}
\newcommand{\bth}[1]{\def\name{Theorem}
\begin{sub}\label{t:#1}}
\newcommand{\blemma}[1]{\def\name{Lemma}
\begin{sub}\label{l:#1}}
\newcommand{\bcor}[1]{\def\name{Corollary}
\begin{sub}\label{c:#1}}
\newcommand{\bdef}[1]{\def\name{Definition}
\begin{sub}\label{d:#1}}
\newcommand{\bprop}[1]{\def\name{Proposition}
\begin{sub}\label{p:#1}}

\newcommand{\rth}[1]{Theorem~\ref{t:#1}}
\newcommand{\rlemma}[1]{Lemma~\ref{l:#1}}

\newcommand{\rprop}[1]{Proposition~\ref{p:#1}}
\newcommand{\BA}{\begin{array}}
\newcommand{\EA}{\end{array}}
\newcommand{\BAN}{\renewcommand{\arraystretch}{1.2}
\setlength{\arraycolsep}{2pt}\begin{array}}
\newcommand{\BAV}[2]{\renewcommand{\arraystretch}{#1}
\setlength{\arraycolsep}{#2}\begin{array}}
\newcommand{\BSA}{\begin{subarray}}
\newcommand{\ESA}{\end{subarray}}
\newcommand{\BAL}{\begin{aligned}}
\newcommand{\EAL}{\end{aligned}}
\newcommand{\BALG}{\begin{alignat}}
\newcommand{\EALG}{\end{alignat}}
\newcommand{\BALGN}{\begin{alignat*}}
\newcommand{\EALGN}{\end{alignat*}}
\newcommand{\note}[1]{\textit{#1.}\hspace{2mm}}
\newcommand{\Proof}{\note{Proof}}
\newcommand{\qeda}{\hspace{10mm}\hfill $\square$}
\newcommand{\Remark}{\note{Remark}}

\newcommand{\norm}[1]{\left \|#1\right \|}

\def\angb<#1>{\langle #1 \rangle}

\newcommand{\myfrac}[2]{{\displaystyle \frac{#1}{#2} }}
\newcommand{\myint}[2]{{\displaystyle \int_{#1}^{#2}}}

\newcommand{\prt}{\partial}

\newcommand{\ti}{\times}

\newcommand{\nind}{\noindent}

\def\ga{\alpha}     \def\gb{\beta}       
       \def\gd{\delta}      \def\ge{\epsilon}
\def\gth{\theta}                         
\def\gf{\phi}           
            \def\gl{\lambda}
\def\gm{\mu}        \def\gn{\nu}         \def\gp{\pi}
            
\def\gs{\sigma}       \def\gt{\tau}
      \def\gw{\omega}
                \def\gz{\zeta}
\def\Gg{\Gamma}     \def\Gd{\Delta}      
\def\Gth{\Theta}
\def\Gl{\Lambda}

   \def\CE{{\mathcal E}}   
      
      \def\CL{{\mathcal L}}

   \def\BBN {\mathbb N}    
   \def\BBR {\mathbb R}    
       
   \def\BBX {\mathbb X}

\title{Nonlinear boundary value problems \\relative to the one dimensional 
heat equation$\phantom{d^{q^p}}$}
\author{Laurent V\'eron}

\date{Month dd, yyyy\\ Revised Month dd, yyyy\\Accepted Month dd, yyyy}
\dedicated{To Julian with high esteem and sincere friendship}

\intesta

\begin{document}

\maketitle

\keywords{Nonlinear heat flux, Singularities, Radon measures, Marcinkiewicz spaces}{35J65, 35L71}
\abstract {We consider the problem of existence of a solution $u$ to $\prt _tu-\prt _{xx} u=0$ in $(0,T)\times\BBR_+ $ subject to the boundary condition 
$-u_x(t,0)+g(u(t,0))=\gm$ on $(0,T)$ where $\gm$ is a measure on $(0,T)$ and $g$ a continuous nondecreasing function.  When $p>1$ we study the set of self-similar solutions of  $\prt _tu-\prt _{xx} u=0$ in $\BBR_+\ti\BBR_+$ such that  $-u_x(t,0)+u^p=0$ on $(0,\infty)$. At end, we present various extensions to a higher dimensional framework.}

\section{Introduction}
\setcounter{equation}{0}
Let  $g:\BBR\mapsto\BBR$ be a continuous nondecreasing  function. Set $Q^T_{\BBR_+}=(0,T)\ti \BBR_+$ for $0<T\leq\infty$ and  
$\prt_\ell Q^T_{\BBR_+}=\overline{\BBR_+}\ti\{0\}$. 
The aim of this article is to study the following 1-dimensional heat equation with a nonlinear flux on the parabolic boundary 
\bel{I-1-0}\BA {lll}
\phantom{g(u---}
\phantom{,,,,}u_t-u_{xx}=0\qquad&\text{in }\;Q^T_{\BBR_+}\\[1mm]
\phantom{;}-u_x(.,0)+g(u(.,0))=\gm\qquad&\text{in }\;[0,T)\\[1mm]\phantom{-------,,}
u(0,.)=\gn\qquad&\text{in }\;\BBR_+,
\EA
\ee
where $\gn,\gm$ are Radon measures in $\BBR_+$ and  $[0,T)$ respectively.  
A related problem in  $Q^\infty_{\BBR_+}$ for which there exist explicit solutions is the following,
\bel{I-1-2}\BA {lll}
\phantom{-------,,ii}
u_t-u_{xx}=0\quad&\text{in }\;Q^\infty_{\BBR_+}\\[0mm]\phantom{--}
-u_x(t,0)+|u|^{p-1}u(t,0)=0&\text{for all }\;t>0\\[1mm]
\phantom{-------,,i}\displaystyle
\lim_{t\to 0}u(t,x)=0&\text{for all }\;x>0,
\EA
\ee
 where $p>1$. Problem $(\ref{I-1-2})$ is invariant under the transformation $T_k$ defined for all $k>0$ by 
 \bel{I-1-3}\BA {lll}
T_k[u](t,x)=k^{\frac{1}{p-1}}u(k^2t,kx).
\EA
\ee
This leads naturaly to look for existence of self-similar solutions under the form
 \bel{I-1-4}\BA {lll}
u_s(t,x)=t^{-\frac{1}{2(p-1)}}\gw\left(\frac x{\sqrt t}\right).
\EA
\ee
Putting $\eta=\frac x{\sqrt t}$, $\gw$ satisfies 
 \bel{I-1-5}\BA {lll}
\!\!-\gw''-\myfrac12\eta \gw'-\myfrac1{2(p-1)}\gw=0\quad&\text{in }\; \BBR_+\\[2mm]
\phantom{--}
-\gw'(0)+|\gw|^{p-1}\gw(0)=0\\[1mm]
\phantom{-u(tx),,-}
\displaystyle
\lim_{\eta\to \infty}\eta^{\frac{1}{p-1}}\gw(\eta)=0.
\EA
\ee 
Self-similar solutions of non-linear diffusion equations such as porous-media or fast-diffusion equation were discovered long time ago by Kompaneets and Zeldovich and a thourougful study was made by Barenblatt, reducing the study to the one of integrable ordinary differential equations with explicit solutions. Concerning semilinear heat equation Brezis, Terman and Peletier opened the study of self-similar solutions of semilinear heat equations in proving in \cite{BPT}  the existence of a positive strongly singular function satisfying 
 \bel{I-1-6}\BA {lll}
u_t-\Gd u+|u|^{p-1}u=0\quad&\text{in }\;\BBR_+\ti\BBR^n,
\EA
\ee 
and vanishing at $t=0$ on $\BBR^n\setminus\{0\}$. They called it the  {\it very singular solution}. Their method of construction is based upon the study of an ordinary differential equation with a phase space analysis. A new and more flexible method based upon variational analysis has been provided by \cite{EK}. Other singular solutions of 
$(\ref{I-1-6})$ in  different configurations such as boundary singularities have been studied in \cite{MV1}. We set $K(\eta)=e^{\eta^2/4}$ and 

 \bel{I-1-6'}\BA {lll}
L^2_K(\BBR_+)=\left\{\phi\in L^1_{loc}(\BBR_+):\myint{\BBR_+}{}\gf^2 Kdx:=\norm\gf^2_{L^2_K}<\infty\right\},
\EA
\ee 
and, for $k\geq1$,
 \bel{I-1-6''}\BA {lll}\displaystyle
H^{k}_K(\BBR_+)=\left\{\phi\in L^2_K(\BBR_+):\sum_{\ga=0}^k\norm{\gf^{(\ga)}}^2_{L^2_K}:=\norm\gf^2_{H^{k}_K}<\infty\right\}.
\EA
\ee
Let us denote by $\CE$ the subset  of $H^{1}_K(\BBR_+)$ of weak solutions of $(\ref{I-1-5})$ that is the set of functions satisfying
 \bel{I-1-6'''}\BA {lll}\displaystyle
\myint{0}{\infty}\left(\gw'\gz'-\myfrac{1}{2(p-1)}\gw \gz\right) K(\eta)d\eta+\left(|\gw|^{p-1}\gw\gz\right)(0) =0,
\EA
\ee
and by $\CE_+$ the subset of nonnegative solutions. The next result gives the structure of $\CE$.

\bth{vssth} 1- If $p\geq2$, then $\CE=\{0\}$.\smallskip

\noindent 2- If $1<p\leq \frac{3}{2}$, then $\CE_+=\{0\}$\smallskip

\noindent 3 - If $\frac{3}{2}<p<2$ then $\CE=\{\gw_s,-\gw_s,0\}$ where $\gw_s$ is the unique positive solution of $(\ref{I-1-5})$. Furthermore 
 there exists 
$c>1$ such that 
 \bel{I-1-7}\BA {lll}
c^{-1}\eta^{\frac{1}{p-1}-1}\leq e^{\frac{\eta^2}{4}}\gw_s(\eta)\leq c\eta^{\frac{1}{p-1}-1}\;\text{ for all }\;\eta>0.
\EA
\ee
\es
Whenever it exists the function $u_s$ defined in  $(\ref{I-1-4})$ is the limit, when $\ell\to\infty$ of the positive solutions $u_{\ell\gd_0}$  of 
\bel{I-1-8}\BA {lll}
\phantom{--g(u-----)}
\!u_t-u_{xx}=0\quad&\text{in }\;Q^\infty_{\BBR_+}\\[1mm]\phantom{,--}
\!-u_{x}(t,.)+|u|^{p-1}u(t,.)=\ell\gd_0\quad&\text{in }\;[0,T)\\[1mm]
\phantom{--------;}\displaystyle
\lim_{t\to 0}u(t,x)=0&\text{for all }\;x\in\BBR_+.
\EA
\ee
When such a function $u_s$ does not exits the sequence $\{u_{\ell\gd_0}\}$ tends to infinity. This is a charateristic phenomenon of an underlying fractional diffusion associated to the linear equation 
 \bel{I-1-8+}\BA {lll}
\phantom{g(u---}
\phantom{,,,,}u_t-u_{xx}=0\qquad&\text{in }\;Q^\infty_{\BBR_+}\\[1mm]
\phantom{,,+g(u(.,0))}-u_x(.,0)=\gm\qquad&\text{in }\;[0,\infty)\\[1mm]\phantom{-------,,}
u(0,.)=0\qquad&\text{in }\;\BBR_+.
\EA
\ee

More generaly we consider problem $(\ref{I-1-0})$. We define the set $\BBX(Q^T_{\BBR_+})$ of test functions by
 \bel{I-1-9}\BA {lll}
\BBX(Q^T_{\BBR_+})=\left\{\gz\in C_c^{1,2}([0,T)\ti[0,\infty)):\gz _{x}(t,0)=0\;\text{ for }t\in [0,T]\right\}.
\EA
\ee

\medskip

\bdef{weak} Let $\gn,\gm$ be Radon measures in $\BBR_+$ and  $[0,T)$ respectively.  A function $u$ defined in  $\overline{Q^T_{\BBR_+}}$ and belonging to $L^1_{loc}(\overline{Q^T_{\BBR_+}})\cap L^1(\prt_\ell Q^T_{\BBR_+};dt)$
 such that $g(u)\in L^1(\prt_\ell Q^T_{\BBR_+};dt)$ is a weak solution of $(\ref{I-1-0})$ if for every $\gz\in\BBX(Q^T_{\BBR_+})$ there holds
 \bel{I-1-10}\BA {lll}
-\myint{0}{T}\myint{0}{\infty}(\gz_t+\gz_{xx})udxdt+\myint{0}{T}\left(g(u)\gz\right)(t,0)dt\\[4mm]
\phantom{---\myint{0}{T}\myint{0}{\infty}(\gz_t+\gz_{xx})udxdt}=\myint{0}{\infty}\gz d\gn(x)+
\myint{0}{T}\gz (t,0) d\gm(t).
\EA
\ee
\es
We denote by $E(t,x)$ the Gaussian kernel in $\BBR_+\ti\BBR$. The solution of  
\bel{I-1-11}\BA {lll}
v_t- v_{xx}=0\qquad&\text{in }\;Q^\infty_{\BBR_+}\\[1mm]\phantom{v_t-}
-v_{x}=\gd_0\qquad&\text{in }\;\overline{\BBR_+}\\[1mm]\phantom{-}
\phantom{,,}v(0,.)=0\qquad&\text{in }\;\BBR_+,
\EA
\ee
has explicit expression
\bel{I-1-12}\BA {lll}
v(t,x)=2E(t,x)=\myfrac{1}{\sqrt{\gp t}}e^{-\frac{x^2}{4t}}.
\EA
\ee
If $x,y>0$ and $s<t$ we set  $\tilde E(t-s,x,y)=E(t-s,x-y)+E(t-s,x+y)$.  When $\gn\in \mathfrak M^b(\BBR_+)$ and $\gm\in \mathfrak M^b(\overline{\BBR_+})$ the solution of 
\bel{I-1-13}\BA {lll}
\phantom{i-(t,0)}v_t-v_{xx} =0\qquad&\text{in }\;Q^\infty_{\BBR_+}\\[1mm]\phantom{-\prt_{xx} E}
-v_{x}(.,0)=\gm\qquad&\text{in }\;\overline{\BBR_+}\\[1mm]\phantom{--}
\!\phantom{i-(t,0)}u(0,.)=\gn\qquad&\text{in }\;\BBR_+,
\EA
\ee
is given by
\bel{I-1-14}\BA {lll}
v_{\gn,\gm}(t,x)=\myint{0}{\infty}\tilde E(t,x,y)d\gn(y)+2\myint{0}{t}E(t-s,x) d\gm(s)\\[3mm]
\phantom{v_{\gn,\gm}(t,x)}=\CE_{\BBR_+}[\gn](t,x)+\CE_{\BBR_+\ti\{0\}}[\gm](t,x)=\CE_{Q^\infty_{\BBR_+}}[(\gn,\gm)](t,x).
\EA
\ee

We prove the following existence and uniqueness result.
\bth{exist-uniq} Let  $g:\BBR\mapsto\BBR$ be a continuous nondecreasing function such that $g(0)=0$.  If $g$ satisfies 
 \bel{I-1-15}\BA {lll}
\myint{1}{\infty}(g(s)-g(-s))s^{-3}ds<\infty,
\EA
\ee
then for any bounded Borel measures $\gn$ in $\BBR_+$ and $\gm$ in $[0,T)$, there exists a unique weak solution $u:=u_{\gn,\gm}\in L^1({Q^T_{\BBR_+}})$ of $(\ref{I-1-0})$. Furthermore the mapping $(\gn,\gm)\mapsto u_{\gn,\gm}$ is nondecreasing.
\es

When $g(s)=|s|^{p-1}s$, condition $(\ref{I-1-15})$ is satisfied if 
 \bel{I-1-16}\BA {lll}
0<p<2.
\EA
\ee

The above result is still valid under minor modifications if $\BBR_+$ is replaced by a bounded interval $I:=(a,b)$, and  problem $(\ref{I-1-0})$ by 
\bel{I-1-17}\BA {lll}
\phantom{g(u(}
\phantom{,,t,1---)}u_t-u_{xx}=0\qquad&\text{in }\;Q^T_{I}\\[1mm]\phantom{-,--}
u_x(.,b)+g(u(.,b))=\gm_1\qquad&\text{in }\;[0,T)\\[1mm]\phantom{,}
-u_x(.,a)+g(u(.,a))=\gm_2\qquad&\text{in }\;[0,T)
\\[1mm]\phantom{---------}
u(0,.)=\gn\qquad&\text{in }\;(a,b),
\EA
\ee
where $\gn,\gm_j$ ($j=1,2$) are Radon measures in $I$ and  $(0,T)$ respectively.  \smallskip

In the last section we present the scheme of the natural extensions of this problem to a multidimensional framework 
 \bel{I-1-18}\BA {lll}
\phantom{g(u;}
\phantom{,--,}u_t-\Gd u=0\qquad&\text{in }\;Q^T_{\BBR^n_+}\\[1mm]\phantom{,,,u}
-u_{x_n}+g(u)=\gm\qquad&\text{in }\;\prt_\ell Q^T_{\BBR^n_+}\\[1mm]\phantom{--,,---}
u(0,.)=\gn\qquad&\text{in }\;\BBR^n_+,
\EA
\ee
The construction of solutions with measure data can be generalized but there are some difficulties in the obtention of self-similar solutions. The equation with a source flux
 \bel{I-1-19}\BA {lll}
\phantom{g(u;}
\phantom{,--,}u_t-\Gd u=0\qquad&\text{in }\;Q^T_{\BBR^n_+}\\[1mm]\phantom{,,,,-u}
u_{x_n}+g(u)=0\qquad&\text{in }\;\prt_\ell Q^T_{\BBR^n_+}\\[1mm]\phantom{--,,---}
u(0,.)=\gn\qquad&\text{in }\;\BBR^n_+,
\EA
\ee
has been studied by several authors, in particular Fila, Ishige, Kawakami and Sato \cite{Fil},  \cite{Ish-Ka}, \cite{Ish-Sa}. Their main concern deals with global existence of solutions. \medskip

\nind {\it Aknowledgements}. The author is grateful to the reviewer for mentioning reference \cite{FQ} which pointed out the role of Whittaker's equation which was  used for analyzing the blow-up of positive solutions of $(\ref{I-1-19})$ when $g(u)=u^p$ when $n=1$.
\section{Self-similar solutions}
\subsection{The symmetrization}
We define the operator $\CL_K$ in $C^2_0(\BBR)$ by 
$$\CL_K(\phi)=-K^{-1}(K\phi')'.
$$
The operator $\CL_K$ has been thouroughly studied in \cite{EK}. In particular
\bel{2-1-3}\inf\left\{\myint{-\infty}{\infty}\phi'^2K(\eta)\eta:\myint{-\infty}{\infty}\phi^2K(\eta)d\eta=1\right\}=\myfrac{1}{2}.
\ee
The above infimum is achieved by $\phi_1=(4\gp)^{-\frac 12}K^{-1}$ and 
$\CL_K$ is an isomorphism from $H^1_K(\BBR)$ onto its dual $(H^1_K(\BBR))'\sim H^{-1}_K(\BBR)$. Finally $\CL_K^{-1}$ is compact from 
 $L^2_K(\BBR)$ into $H^1_K(\BBR)$, which implies that $\CL_K$ is a Fredholm self-adjoint operator with
 $$\gs(\CL_K)=\left\{\gl_j=\tfrac{1+j-1}{2}:j=1, 2,...\right\},$$ 
 and 
 $$ker\left(\CL_K-\gl_jI_d\right)=span\left\{\phi_1^{(j)}\right\}.
 $$
If $\gf$ is defined in $\BBR_+$, $\tilde \gf(x)=\gf(-x)$ is the symmetric with respect to $0$ while $\phi^*(x)=-\phi(-x)$ is the antisymmetric with respect to $0$. The operator $\CL_K$ restricted to $\BBR_+$ is denoted by $\CL^+_K$. The operator $\CL^{+,N}_{K}$ with Neumann condition at $x=0$ is again a Fredholm operator. This is also valid for the operator $\CL^{+,D}_{K}$ with Dirichlet condition at $x=0$. Hence, if $\phi$ is an eigenfunction of $\CL^{+,N}_{K}$, then $\tilde\phi$ is an eigenfunction of $\CL_K$ in $L^2_K(\BBR)$. Similarly, if $\phi$ is an eigenfunction of $\CL^{+,D}_{K}$, then $\phi^*$ is an eigenfunction of $\CL_K$ in $L^2_K(\BBR)$. Conversely, any even (resp. odd) eigenfunction of $\CL_K$  in $L^2_K(\BBR)$ satisfies Neumann (resp. Dirichlet) boundary condition at $x=0$. Hence its restiction to 
$L^2_K(\BBR_+)$ is an eigenfunction of $\CL^{+,N}_{K}$ (resp. $\CL^{+,D}_{K}$). Since $\phi_1^{(j)}$ is even (resp. odd) if and only if 
$j$ is even (resp. odd), we derive
\bel{2-1-1}
H^{1,0}_K(\BBR_+)=\bigoplus_{\ell=1}^\infty span\left\{\phi_1^{(2\ell+1)}\right\},
\ee
and
\bel{2-1-2}
H^{1}_K(\BBR_+)=\bigoplus_{\ell=0}^\infty span\left\{\phi_1^{(2\ell)}\right\}.
\ee

Note that $\phi\in H^1_K(\BBR_+)$ such that $\phi_x(0)=0$ (resp. $\phi(0)=0$) implies $\tilde\phi\in H^1_K(\BBR)$ (resp. $\phi^*\in H^1_K(\BBR)$).
 Furthermore, $\phi_1$ is an eigenfunction of $\CL^+_K$ in $H^1_K(\BBR_+^n)$ with Neumann boundary condition on $\prt\BBR_+^n$ while $\prt_{x_n}\phi_1$ is an eigenfunction of $\CL^+_K$ in $H^1_K(\BBR_+^n)$ with Dirichlet boundary condition on $\prt\BBR_+^n$. We list below two important properties of $H^1_K(\BBR_+)$ valid for any $\gb>0$. Actually they are proved in \cite[Prop. 1.12]{EK} with $H^1_{K^\gb}(\BBR)$ but the proof is valid with $H^1_{K^\gb}(\BBR_+)$.
\bel{imbed}\BA{lll}
\!\!\!\!\!\!(i) \quad\, &\phi\in H^1_{K^\gb}(\BBR_+)\Longrightarrow  K^{\frac{\gb}{2}}\phi\in C^{0,\frac{1}{2}}(\BBR_+)\\[1mm]
\!\!\!\!\!\!(ii)\qquad &H^1_{K^\gb}(\BBR_+)\hookrightarrow L^2_{K^\gb}(\BBR_+)\quad\text{is compact for all }n\geq 1.
\EA
\ee

\noindent


\subsection{Proof of \rth{vssth}-(i)-(ii)}
Assume $p\geq 2$, then $\frac{1}{2(p-1)}\leq \frac 12$. If $\gw$ is a weak solution, then 
$$
\myint{0}{\infty}\left(\gw'^2-\frac{1}{2(p-1)}\gw^2\right)Kd\eta+|\gw|^{p+1}(0)= 0.
$$ 
If  $\frac 12> \frac{1}{2(p-1)}$ we deduce that $\gw=0$. Furthermore, when $\frac 12= \frac{1}{2(p-1)}$ then 
$$|\gw|^{p+1}(0)=  0.
$$
If $\gw$ is nonzero, it is an eigenfunction of $\CL^{+,D}_K$. Since the first eigenvalue is $1$ it would imply $1=\frac{1}{2(p-1)}\leq \frac 12$, contradiction. 
\smallskip

\noindent Assume $1<p\leq \frac3{2}$ and $\gw$ is a nonnegative weak solution. 
We take $\gz(\eta)=\eta e^{-\frac{\eta^2}{4}}=-2\phi'_{1\,}(\eta)$, then  
$$\BA {lll}
\myint{0}{\infty}\left(-\gz''-\myfrac{1}{2(p-1)}\gz\right)\gw  K(\eta)d\eta+\gz'(0)\gw^p(0) =0.
\EA$$
Since $-\gz''=\gz\lfloor_{\BBR_+}>0$ and $\gz'(0)=\phi_1(0)=1$, we derive $\gw\gz=0$ if $1>\frac{1}{2(p-1)}$ and $\gw(0)=0$ if 
$1=\frac{1}{2(p-1)}$. Hence $\gw'(0)=0$ by the equation and $\gw\equiv 0$ by the Cauchy-Lipschitz theorem. \qeda

\subsection{Proof of \rth{vssth}-(iii)}

 We define the following functional on $H^1_K(\BBR^n_+)$
   \bel{2-2-1}
   J(\phi)=\myfrac{1}{2}\myint{0}{\infty}\left(\gf'^2-\myfrac{1}{2(p-1)}\gf^2\right) Kd\eta+\myfrac{1}{p+1}|\gf(0)|^{p+1}.
\ee

\blemma{funct} The functional $J$ is lower semicontinuous in $H^1_{K}(\BBR_+)$. It tends to infinity at infinity and achieves negative values.
\es
\Proof
We write
$$J(\psi)=J_1(\psi)-J_2(\psi)=J_1(\psi)-\frac{1}{2(p-1)}\norm\psi_{L^2_K}^2.$$ 
Clearly $J_1$ is convex and $J_2$ is continuous in the weak topology of $H^1_K(\BBR_+)$ since the imbedding of 
$H^1_K(\BBR_+)$ into $L^2_K(\BBR_+)$ is compact. Hence $J$ is weakly semicontinuous in  $H^1_K(\BBR_+)$.
\smallskip

\noindent Let $\ge>0$, then 
$$J(\ge\phi_1)=\left(\myfrac{1}{4}-\myfrac{1}{4(p-1)}\right)\myfrac{\ge^2\sqrt\gp}{2}+\myfrac{\ge^{p+1}}{p+1}.
$$
Since $1<p<2$, $\frac{1}{4}-\frac{1}{4(p-1)}<0$. Hence $J(\ge\phi_1)<0$ for $\ge$ small enough, thus $J$ achieves negative values on $H^1_{K}(\BBR_+)$. \smallskip

\noindent If $\psi\in H^1_{K}(\BBR_+)$ it can be written in a unique way under the form  $\psi=a\phi_1+\psi_1$ where $a=2\sqrt\gp\psi(0)$ and $\psi_1\in H^{1,0}_K(\BBR_+)$.
Hence, for any $\ge>0$, 
$$\BA {lll}J(\psi)=\myfrac{1}{2}\myint{0}{\infty}\left(\psi_1'^2-\myfrac{1}{2(p-1)}\psi_1^2\right) Kd\eta
+\myfrac{a^2}{2}\myint{0}{\infty}\left(\phi_1'^2-\myfrac{1}{2(p-1)}\phi_1^2\right) Kd\eta
\\[4mm]
\phantom{J(\psi)=}
+a\myint{0}{\infty}\left(\psi_1'\phi_1'-\myfrac{1}{2(p-1)}\psi_1\phi_1\right) Kd\eta
+\myfrac{1}{p+1}|a|^{p+1}
\\[4mm]\phantom{J(\psi)}
\geq \myfrac{2p-3}{4(p-1)}\myint{0}{\infty}\psi_1'^2 Kd\eta-\myfrac{a\ge}{2}\myint{0}{\infty}\left(\psi_1'^2+\myfrac{1}{2(p-1)}\psi_1^2\right) K d\eta\\[4mm]\phantom{J(\psi)=}
+\myfrac{a^2(p-2)\sqrt\gp}{4(p-1)}-\myfrac{ap\sqrt\gp}{4(p-1)\ge}+\myfrac{1}{p+1}|a|^{p+1}.
\EA$$
Note that $\norm{\psi}^2_{H^1_K}\leq 4\left(\norm{\psi'_1}^2_{L^2_K}+a^2\right)$. Since $2p-3>0$, we can take $\ge>0$ small enough in order that
 \bel{2-2-9}\displaystyle
 \lim_{\norm\psi_{H^1_K}\to\infty}J(\psi)=\infty. 
\ee
\qeda\medskip

\noindent{\it End of the proof of \rth{vssth}-(iii)}.
By \rlemma{funct} the functional $J$ achieves its minimum in $H^{1}_K(\BBR_+)$ at  some $\gw_s\neq 0$, and $\gw_s$ can be assumed to be nonnegative since $J$ is even. By the strong maximum principle $\gw_s>0$, and by the method used in the proof of \cite[Proposition 1] {MoV} is is easy to prove that positive solutions belong to $H^2_K(\BBR_+)$. Assume that $\tilde \gw$ is another positive solution, then
$$\myint{0}{\infty}\left(\myfrac{(K\gw'_s)'}{\gw_s}-\myfrac{(K\tilde\gw'_s)'}{\tilde\gw_s}\right) (\gw_s^2-\tilde\gw_s^2)d\eta=0.
$$
Integration by parts, easily justified by regularity, yields 
$$\BA {lll}
\myint{0}{\infty}\left(\myfrac{(K\gw'_s)'}{\gw_s}-\myfrac{(K\tilde\gw'_s)'}{\tilde\gw_s}\right) (\gw_s^2-\tilde\gw_s^2)d\eta\\[4mm]
\phantom{-----}
=\left[K\gw_s'\left(\gw_s-\myfrac{\tilde\gw_s^2}{\gw_s}\right)-K\tilde\gw_s'\left(\myfrac{\gw_s^2}{\tilde\gw_s}-\tilde\gw_s\right)\right]^\infty_0\\[4mm]
\phantom{-------}-\myint{0}{\infty}\left(\gw_s-\myfrac{\tilde\gw_s^2}{\gw_s}\right)'K\gw_s'd\eta+
\myint{0}{\infty}\left(\myfrac{\gw_s^2}{\tilde\gw_s}-\tilde\gw_s\right)'K\gw_s'd\eta\\[4mm]
\phantom{-----}
=-\left(\gw_s^{p-1}-\tilde\gw_s^{p-1}\right)\left(\gw_s^{2}-\tilde\gw_s^{2}\right)(0)
\\[4mm]
\phantom{-------}
-\myint{0}{\infty}\left(\left(\myfrac{\gw_s'\tilde\gw_s-\gw_s\tilde\gw_s'}{\tilde\gw_s}\right)^2+\left(\myfrac{\gw_s\tilde\gw_s'-\tilde\gw_s\gw_s'}{\gw_s}\right)^2\right)d\eta.
\EA
$$
This implies that $\gw_s=\tilde\gw_s$. The proof of $(\ref{I-1-7})$ is similar as the proof of estimate (2.5) in \cite[Theorem 4.1]{MV1}.\qeda 

\subsection{The explicit approach}
This part is an adaptation to our problem of what has been done in \cite{FQ} concerning the blow-up problem in equation $(\ref{I-1-19})$. 
Let $\gw$ be a solution of 
\bel{W-0}
\gw''+\myfrac12\eta \gw'+\myfrac1{2(p-1)}\gw=0\qquad\text{in }\; \BBR_+.
\ee
We set 
$$r=\frac{\eta^2}{4} \, \text{ and }\; \gw(\eta)=r^{-\frac{1}{4}}e^{-\frac{r}{2}}Z(r).
$$          
Then $Z$ satisfies the Whittaker equation (with the standard notations)
\bel{W-1}
Z_{rr}+\left(-\myfrac{1}{4}+\myfrac{k}{r}+\myfrac{1-4\gm^2}{4r^2}\right)Z=0
\ee
where $k=\frac{1}{2(p-1)}-\frac 14$ and $\gm=\frac 14$. Notice that the only difference with the expression in \cite[Lemma 3.1]{FQ} is 
the value of the coefficient $k$. This equation admits two linearly independent solutions
$$Z_1(r)=e^{-\frac r2}r^{\frac 12+\gm}U\left(\tfrac{1}{2}+\gm-k,1+2\gm,r\right),
$$
and
$$Z_2(r)=e^{-\frac r2}r^{\frac 12+\gm}M\left(\tfrac{1}{2}+\gm-k,1+2\gm,r\right).
$$
The functions $U$ and $M$ are the Whittaker functions which play an important role not only in analysis but also in group theory. The have the following 
asymptotic expansion as $r\to\infty$ (see e.g. \cite{AS}),
$$U\left(\tfrac{1}{2}+\gm-k,1+2\gm,r\right)=r^{k-\gm-\frac 12}\left(1+O(r^{-1}\right)=r^{\frac {1}{2(p-1)}-1}\left(1+O(r^{-1}\right),
$$
and 
$$\BA {lll}
M\left(\tfrac{1}{2}+\gm-k,1+2\gm,r\right)=\myfrac{\Gg(1+2\gm)}{\Gg(\tfrac{1}{2}+\gm-k)}e^rr^{-(\gm+\frac12+k)}\left(1+O(r^{-1}\right)\\[4mm]
\phantom{M\left(\tfrac{1}{2}+\gm-k,1+2\gm,r\right)}
=\myfrac{\Gg(\tfrac 32)}{\Gg(1-\tfrac{1}{2(p-1)})}e^rr^{-\frac p{2(p-1)}}\left(1+O(r^{-1}\right).
\EA$$
Then 
$$Z_1(r)=r^{\frac{1}{2(p-1)}-\frac 14}e^{-\frac{r}{2}}\left(1+O(r^{-1}\right),
$$
and 
$$Z_2(r)=\myfrac{\Gg(\tfrac 32)}{\Gg(1-\tfrac{1}{2(p-1)})}r^{\frac 14-\frac{1}{2(p-1)}-}e^\frac r2\left(1+O(r^{-1}\right).
$$
To this corresponds the two linearly independent solutions $\gw_1$ and $\gw_2$ of $(\ref{W-0})$ with the following behaviour as $\eta\to\infty$,
\bel{W-2}\BA{lll}
(i)\qquad&\gw_1(\eta)=c_1\eta^{\frac{1}{p-1}-1}e^{-\frac{\eta^2}{4}}\left(1+O(\eta^{-2}\right),\qquad\qquad\qquad\\[4mm]

(ii)\qquad&\gw_2(\eta)=c_2\eta^{-\frac{1}{p-1}}\left(1+O(\eta^{-2}\right).
\EA\ee
Clearly only $\gw_1$ satisfies the decay estimate $\gw(\eta)=o(\eta^{-\frac{1}{p-1}})$ as $\eta\to\infty$. Hence the solution $\gw$ is a multiple of 
$\gw_1$ and the multiplicative constant $c$ is adjusted in order to fit the condition $\gw'(0)=\gw^p(0)$.

\section{Problem with measure data}
\subsection{The regular problem}
Set $G(r)=\int_0^rg(s)ds$. We consider the functional $J$ in $L^2(\BBR_+)$ with domain $D(J)=H^1(\BBR_+)$ defined by 
$$ J(u)=\myfrac12\myint{0}{\infty}u_x^2dx+G(v(0)).
$$
It is convex and lower semicontinuous in $L^2(\BBR_+)$ and its subdifferential $\prt J$ sastisfies
$$\myint{0}{\infty}\prt J(u)\gz dx=\myint{0}{\infty}u _x\gz _x dx+g(u(0))\gz(0),
$$
for all $\gz\in H^1(\BBR_+)$. Therefore
$$\myint{0}{\infty}\prt J(u)\gz dx=-\myint{0}{\infty}u_{xx}\gz dx+(g(u(0))-u_x(0))\gz(0).
$$
Hence 
\bel{II-1-1}
\prt J(u)=-u_{xx}\,\;\text{for all }\, u\in D(\prt J)=\{v\in H^1(\BBR_+):v_x(0)=g(v(0))\}.
\ee
The operator $\prt J$ is maximal monotone, hence it generates a semi-group of contractions. Furthermore, for any $u_0\in L^2(\BBR_+)$ and 
$F\in L^2(0,T;L^2(L^2(\BBR_+))$ there exists a unique strong solution to 
\bel{II-1-2}\BA {lll}
U_t+\prt J(U)=F\quad\text{a.e. on }\, (0,T)\\
\phantom{-,,--}U(0)=u_0.
\EA\ee
\bprop{Brcase} Let  $\gm\in H^1(0,T)$ and $\gn\in L^2(\BBR_+)$. Then there exists a unique function 
$u\in C([0,T];L^2(\BBR_+)$ such that $\sqrt tu_{xx}\in L^2((0,T)\ti\BBR_+)$ which satisfies 
$(\ref{II-1-3})$. The mapping $(\gm,\gn)\mapsto u:=u_{\gm,\gn}$ is non-decreasing and $u$ is a weak solution in the sense that it satisfies 
$(\ref{I-1-10})$.
\es
\Proof Let $\eta\in C^2_0([0,\infty))$ such that $\eta (0)=0$, $\eta' (0)=1$. If $f\in H^1(0,T)$, $\gn\in L^2(\BBR_+)$, and $u$ is a solution of 
\bel{II-1-3}\BA {lll}
\phantom{g(u---}
\phantom{,,-}u_t-u_{xx}=0\qquad&\text{in }\;Q^T_{\BBR_+}\\[1mm]
\phantom{,,}\!-u_x(.,0)+g(u(.,0))=\gm(t)\qquad&\text{in }\;[0,T)\\[1mm]\phantom{--------}
u(0,.)=\gn\qquad&\text{in }\;\BBR_+,
\EA
\ee
where $\gn\in L^2(\BBR_+)$, then the function $v(t,x)=u(t,x)-\gm(t)\eta(x)$ satisfies 
\bel{II-1-4}\BA {lll}
\phantom{g(u---}
\phantom{,-,}v_t-v_{xx}=F\qquad&\text{in }\;Q^T_{\BBR_+}\\[1mm]
\phantom{,}-v_x(.,0)+g(v(.,0))=0\qquad&\text{in }\;[0,T)\\[1mm]\phantom{------,,,,}
v(0,.)=\gn-\gm(0)\eta\qquad&\text{in }\;\BBR_+,
\EA
\ee
with $F(t,x)=-(\gm'(t)\eta(x)+\gm(t)\eta''(x))$. The proof of the existence follows by using \cite[Theorem 3.6]{Br-OMM}. \\
\noindent Next, let $(\tilde \gm,\tilde\gn)\in H^1(0,T)\ti L^2(\BBR_+)$ such that $\tilde \gm\leq \gm$ and $\tilde\gn\leq\gn$ and let  $\tilde u=u_{\tilde \gm,\tilde\gn}$, then 
$$\BA {lll}\myfrac 1{2}\myfrac d{dt}\myint{0}{\infty}(\tilde u-u)^2_+dx+\myint{0}{\infty}\left(\prt_x(\tilde u-u)_+\right)^2dx-\left(\tilde \gm(t)-\gm(t)\right)(\tilde u(t,0)-u(t,0))_+\\[4mm]\phantom{-----------}
+\left(g(\tilde u(t,0))-g(u(t,0))\right))(\tilde u(t,0)-u(t,0))=0.
\EA$$
Then 
$$\myint{0}{\infty}(\tilde u-u)^2_+dx\lfloor_{t=0}\;\Longrightarrow\myint{0}{\infty}(\tilde u-u)^2_+dx=0\quad\text{ on }\,[0,T].$$ 
We can also use $(\ref{I-1-14})$ to express the solution of $(\ref{II-1-3})$:
$$u(t,x)=\myint{0}{\infty}\tilde E(t,x,y)\gn(y)dy+2\myint{0}{t}E(t-s,x) (\gm(s)-g(u(s,0)))ds.
$$
In particular, if $g(0)=0$, then 
$$|u(t,x)|\leq \myint{0}{\infty}\tilde E(t,x,y)|\gn(y)|dy+2\myint{0}{t}E(t-s,x) |\gm(s)|ds.
$$
The proof of $(\ref{I-1-10})$ follows  since $u$ is a strong solution. 
\qeda
\medskip

Next, we prove that the problem is well-posed if $\gm\in L^1(0,T)$. 

\bprop{L^1} Assume $\{\gn_n\}\subset C_c(\BBR_+)$ and $\{\gm_n\}\subset C^1([0,T])$ are Cauchy sequences in $L^1(\BBR_+)$ and $L^1(0,T)$ respectively. Then the sequence $\{u_n\}$ of solutions of 
\bel{II-1-5}\BA {lll}\phantom{------}
\!u_{n\,t}-u_{n\,xx}=0\qquad&\text{in }\;Q^T_{\BBR_+}\\[1mm]
\,\,-u_{n\,x}(.,0)+g(u_n(.,0))=\gm_n(t)\qquad&\text{in }\;[0,T)\\[1mm]\phantom{--------}
\!u_n(0,.)=\gn_n\qquad&\text{in }\;\BBR_+,
\EA
\ee
converges in $C([0,T];L^1(\BBR_+)$ to a function $u$ which satisfies $(\ref{I-1-10})$.  
\es
\Proof For $\ge>0$ let $p_\ge$ be an odd $C^1$ function defined on $\BBR$ such that $p_\ge'\geq 0$ and $p_\ge(r)=1$ on $[\ge,\infty)$, and put $j_\ge(r)=\int_0^rp_\ge(s)ds$. Then 
$$\BA {lll}\myfrac {d}{dt}\myint{0}{\infty}j_\ge(u_n-u_m)dx+\myint{0}{\infty}(u_{n\,x}-u_{m\,x})^2p'_\ge(u_n-u_m)dx\\[4mm]\phantom{--------}
+\left(g(u_n(t,0))-g(u_m(t,0))\right)p_\ge(u_n(t,0)-u_m(t,0))\\[4mm]
\phantom{-----------}=\left(\gm_n(t)-\gm_m(t)\right)p_\ge(u_n(t,0)-u_m(t,0)).
\EA$$
Hence
 $$\BA {lll}\myint{0}{\infty}j_\ge(u_n-u_m)(t,x)dx+\left(g(u_n(t,0))-g(u_m(t,0))\right)p_\ge(u_n(t,0)-u_m(t,0)) \\[4mm]
 \phantom{------}
\leq\myint{0}{\infty}j_\ge(\gn_n-\gn_m)dx+\left(\gm_n(t)-\gm_m(t)\right)p_\ge(u_n(t,0)-u_m(t,0)).
 \EA$$
 Letting $\ge\to 0$ implies $p_\ge\to sgn_0$, hence for any $t\in [0,T]$, 
 \bel{II-1-6}\BA {lll}\myint{0}{\infty}|u_n-u_m|(t,x)dx+|g(u_n(t,0))-g(u_m(t,0)| \\[4mm]
 \phantom{-----------}
\leq\myint{0}{\infty}|\gn_n-\gn_m|dx+|\gm_n(t)-\gm_m(t)|.
\EA
\ee
Therefore $\{u_n\}$ and $\{g(u_n(.,0)\}$ are Cauchy sequences in $C([0,T];L^1(\BBR_+))$ and  $C([0,T])$ respectively with limit $u$ and $g(u)$
and $u=u_{\gn,\gm}$ satisfies $(\ref{I-1-10})$. 
If we assume that $(\gn,\tilde\gn)$ and $(\gm,\tilde\gm)$ are couples of elements of $L^1(\BBR_+)$ and $L^1(0,T)$ respectively and if 
$u=u_{\gn,\gm}$ and $\tilde u=u_{\tilde\gn,\tilde\gm}$, there holds by the above technique, 
 \bel{II-1-7}\BA {lll}\myint{0}{\infty}|u-\tilde u|(t,x)dx+|g(u(t,0))-g(\tilde u(t,0)| \\[4mm]
 \phantom{------}
\leq\myint{0}{\infty}|\tilde\gn-\tilde\gn|dx+|\tilde\gm(t)-\tilde\gm(t)|\quad\text{for all }\,t\in [0,T].
\EA
\ee
\qeda

The following lemma is a parabolic version of an inequality due to Brezis.
\blemma{Brtype} Let $\gn \in L^1(\BBR_+)$ and $\gm \in L^1(0,T)$ and $v$ be a function  defined in $[0,T)\ti\BBR_+$, belonging to $L^1(Q^T_{\BBR_+})\cap L^1(\prt_\ell Q^T_{\BBR_+})$ and satisfying 
\bel{II-1-8}\BA {lll}
-\myint{0}{T}\!\!\myint{0}{\infty}(\gz_t+\gz_{xx})vdxdt=\myint{0}{T}\gz(.,0) \gm dt+\myint{0}{\infty}\gn\gz dx.
\EA
\ee
Then for any $\gz\in\BBX(Q^T_{\BBR_+})$, $\gz\geq 0$, there holds
\bel{II-1-9}\BA {lll}
-\myint{0}{T}\!\!\myint{0}{\infty}(\gz_t+\gz_{xx})|v|dxdt\leq\myint{0}{\infty}\gz(.,0) sign (v)\gm dt+\myint{0}{\infty}|\gn|\gz dx.
\EA\ee
Similarly
\bel{II-1-10}\BA {lll}
-\myint{0}{T}\!\!\myint{0}{\infty}(\gz_t+\gz_{xx})v_+dxdt\leq\myint{0}{\infty}\gz(.,0) sign_+ (v)\gm dt+\myint{0}{\infty}\gn_+\gz dx.
\EA\ee
\es
\Proof Let $p_\ge$ be the approximation of $sign_0$ used in \rprop{L^1} and $\eta_\ge$ be the solution of 
$$\BA{lll}
-\eta_{\ge\,t}-\eta_{\ge\,xx}=p_\ge(v)\quad&\text{in }Q^T_{\BBR_+}\\[1mm]
\phantom{--,}
\eta_{\ge\,x}(.,0)=0&\text{in }[0,T]\\[1mm]
\phantom{---}
\eta_\ge(0,.)=0&\text{in }\BBR_+.
\EA$$
Then $|\eta_\ge|\leq \eta^*$ where $\eta^*$ satisfies 
$$\BA{lll}
-\eta^*_{t}-\eta^*_{xx}=1\quad&\text{in }Q^T_{\BBR_+}\\[1mm]
\phantom{--}
\eta^*_{x}(.,0)=0&\text{in }[0,T]\\[1mm]
\phantom{--}
\eta^*(0,.)=0&\text{in }\BBR_+.
\EA$$
Although $\eta_\ge$ does not belong to $\BBX(Q^T_{\BBR_+})$ (it is not in $C^{1,2}([0,T)\ti\BBR_+)$, it is an admissible test function and we deduce that there exists a unique solution to $(\ref{II-1-8})$. Thus $v$ is given by expression $(\ref{I-1-14})$.\par

In order to prove $(\ref{II-1-9})$, we can assume that $\gm$ and $\gn$ are smooth, $\gz\in \BBX(Q^T_{\BBR_+})$, $\gz\geq 0$ and set 
$h_\ge=p_\ge(v)\gz$ and $w_\ge=vp_\ge(v)$, then
\bel{II-1-11}\BA {lll}\myint{0}{\infty}h_{\ge\,xx}vdx=\myint{0}{\infty}\left(2p'_\ge(v)v_x\gz_x+p_\ge(v)\gz_{xx}+\gz(p_\ge(v))_{xx}\right)vdx\\[4mm]
\phantom{\myint{0}{\infty}}
=\myint{0}{\infty}\left(2vp'_\ge(v)v_x\gz_x-w_{\ge\, x}\gz_x-(v\gz)_x(p_\ge(v))_x\right)dx\\[0mm]
\phantom{\myint{0}{\infty}------------}
-\gz(t,0)v(t,0)p'_\ge(v(t,0))v_x(t,0)
\\[0mm]
\phantom{\myint{0}{\infty}}
=-\myint{0}{\infty}\left(\gz_x(j_\ge(v))_x+\gz p'(v)_\ge v_x^2\right)dx-\gz(t,0)v(t,0)p'_\ge(v(t,0))v_x(t,0)
\\[0mm]
\phantom{\myint{0}{\infty}}
=-\myint{0}{\infty}\left(\gz p'(v)_\ge v_x^2 -j_\ge(v)\gz_{xx}\right)dx-\gz(t,0)v(t,0)p'_\ge(v(t,0))v_x(t,0),
\EA\ee
and
\bel{II-1-12}\BA {lll}
\myint{0}{T}h_{\ge\,t}vdt=\myint{0}{T}(p_\ge(v)\gz_t+p'_\ge(v)\gz v_t)vdt.
\EA\ee
Since $v$ is smooth
$$\BA {lll}
0=\myint{0}{T}\myint{0}{\infty}(v_t-v_{xx})h_\ge dx dt\\
\phantom{0}
=-\myint{0}{T}\myint{0}{\infty}(h_{\ge\,t}+h_{\ge\,xx})v dx dt-\myint{0}{\infty}h_\ge(0,x)\gn(x)dx\\
\phantom{0---------}
- \myint{0}{T}\left[p_{\ge}(v(t,0))-v(t,0)p'_{\ge}(v(t,0))\right] \gz(t,0)\gm(t)dt.
\EA$$
Therefore, using $(\ref{II-1-9})$ and $(\ref{II-1-10})$,
\bel{II-1-13}\BA {lll}
-\myint{0}{T}\myint{0}{\infty}\left(j_\ge v)\gz_{xx}+vp_{\ge}(v)\gz_t\right)dxdt\\[4mm]
\phantom{---------}+\myint{0}{T}\myint{0}{\infty}\left(\gz p_\ge'(v)v_x^2-vp'_\ge(v)v_t\gz\right) dxdt
\\[4mm]
\phantom{---------}
=\myint{0}{\infty}h_\ge(0,x)\gn(x) dx+\myint{0}{T}h_\ge(t,0)\gm(t) dt.
\EA\ee
Put $\ell_\ge(s)=\int_0^srp'_\ge(r) dr$, then $|\ell_\ge (s)\leq c\ge^{-1}s^2\chi_{_{[-\ge,\ge]}}(s)|$. Since
$$\BA {lll}
\myint{0}{T}\myint{0}{\infty}\gz v p'_\ge(v)v_t dxdt=-\myint{0}{\infty}\ell_\ge(v(0,x))\gz(x)dx-\myint{0}{T}\myint{0}{\infty}\gz_t  \ell_\ge(v)dxdt,
\EA$$
and $\gz$ has compact support, it follows that
$$\displaystyle \lim_{\ge\to 0}\myint{0}{T}\myint{0}{\infty}\gz v p'_\ge(v)v_t dxdt=0.
$$
Letting $\ge\to 0$ in $(\ref{II-1-13})$, we derive $(\ref{II-1-9})$ for smooth $v$. Using \rprop{L^1} completes the proof of $(\ref{II-1-9})$. The proof of $(\ref{II-1-10})$ is similar.\qeda\medskip

\noindent \Remark Inequalities $(\ref{II-1-9})$ and $(\ref{II-1-10})$ hold if $\gz(t,x)$ does not vanish if $|x|\geq R$ for some $R$ but if it satisfies
\bel{II-1-14}\BA {lll}\displaystyle
\lim_{x\to\infty}\sup_{t\in [0,T]}(\gz(t,x)+|\gz_x(t,x)|)=0.
\EA\ee
The proof follows by replacing $\gz(t,x)$ by $\gz(t,x)\eta_n(x)$ where $\eta_n\in C^\infty_c(\BBR_+)$ with $0\leq \eta_n\leq 1$, $\eta_n(x)=1$ on $[0,n]$, 
$\eta_n(x)=0$ on $[n+1,\infty)$, $|\eta'_n|\leq 2$, $|\eta''_n|\leq 4$. Then $\eta_n\gz\in\BBX(Q_{\BBR_+}^T)$ by letting $n\to\infty$ and the proof follows by letting $n\to\infty$.
\subsection{Proof of \rth{exist-uniq}}

We give first some {\it heat-ball} estimates relative to our problem. For $r>0$, $x\in\BBR_+$ and $t\in\BBR$ we set
\bel{III-1-1e}
\BA{lll}
e(t,x;r)=\left\{(s,y)\in (0,T)\ti\BBR_+:s\leq t,\, \tilde E(t-s,x,y)\geq r\right\}.
\EA
\ee
Since 
$$e(t,x;r)\subset [t-\tfrac{1}{4\gp er^2},t]\ti[x-\tfrac{1}{r\sqrt{\gp e}},x+\tfrac{1}{r\sqrt{\gp e}}],$$
there holds 
\bel{III-1-2e}|e(t,x;r)|\leq \myfrac{1}{2r^3(\gp e)^{\frac32}},
\ee
and if 
\bel{III-1-3e}
e^*(t;r)=\left\{s\in (0,T):s\leq t,\,  E(t-s,0,0)\geq r\right\},
\ee
then
we have 
\bel{III-1-4e}e^*(t;r)\subset [t-\tfrac{1}{4\gp er^2},t]\Longrightarrow |e^*(t;r)|\leq \myfrac{1}{4r^2\gp e}.
\ee
If $G$ is a measured space, $\gl$ a positive measure on $G$ and $q>1$, $M^q(G,\gl)$ is the Marcinkiewicz space of measurable functions 
$f:G\mapsto \BBR$ satisfying for some constant $c>0$ and all measurable set $E\subset G$, 
 \bel{IIl-1-4}
 \myint{E}{}|f|d\gl\leq c\left(\gl (E)\right)^{\frac{1}{p'}},
 \ee
 and 
 $$\norm f_{M^q(G,\gl)}=\inf\{c>0\,\text{ s.t. $(\ref{III-1-4e})$ holds}\}.
 $$
\blemma{reg1} Assume $\gm$,$\gn$ are bounded measure in $\overline{\BBR_+}$ and $\BBR_+$ respectively and $u$ is the solution of 
$(\ref{I-1-13})$ given by $(\ref{I-1-14})$ and $v_{\gn,\gm}$ is the solution of $(\ref{I-1-13})$. Then 
\bel{IIl-1-5}
\norm{v_{\gn,\gm}}_{M^3(Q^T_{\BBR_+})}+\norm{v_{\gn,\gm}\lfloor_{\prt Q_{\BBR_+}^T}}_{M^2(\prt Q_{\BBR_+}^T)}\leq c
\left(\norm{\gm}_{\mathfrak M(\prt Q_{\BBR_+}^T)}+\norm\gn_{\mathfrak M(Q^T_{\BBR_+})}\right).
\ee
\es
\Proof 
First we consider $v_{0,\gm}$
$$v_{0,\gm}(t,x)=2\myint{0}{t}E(t-s,x)d\gm(s).$$
If $F\subset [0,T]$ is a Borel set, than for any $\gt>0$
$$\BA {lll}\myint{F}{}E(t-s,0)ds=\myint{F\cap \left\{E\leq \gt\right\}}{}E(t-s,0)ds+\myint{F\cap \left\{E> \gt\right\}}{}E(t-s,0)ds\\[4mm]
\phantom{\myint{F}{}E(t-s,0)ds}\leq \gt|F|+\myint{\left\{E> \gt\right\}}{}E(t-s,0)ds\\[4mm]
\phantom{\myint{F}{}E(t-s,0)ds}\leq \gt|F|-\myint{\gt}{\infty}\gl d|e^*(t,\gl)|\\[4mm]
\phantom{\myint{F}{}E(t-s,0)ds}\leq \gt|F|+\myint{\gt}{\infty}\gl d|e^*(t,\gl)|\\[4mm]
\phantom{\myint{F}{}E(t-s,0)ds}\leq \gt|F|+\myfrac{1}{4\gp e\gt}.
\EA$$
If we choose $\gt^2=\frac{1}{4\gp e|F|}$, we derive
\bel{IIl-1-6}
\myint{F}{}E(t-s,0)ds\leq \myfrac{|F|^{\frac 12}}{\sqrt{\gp e}}.
\ee
If $F\subset (0,T)$ is a Borel set then
$$\left|\myint{F}{}v_{0,\gm}(t,0)dt\right|=2\left|\myint{0}{t}\myint{F}{}E(t-s,0)dtd\gm(s)\right|\leq 
\myfrac{2|F|^{\frac 12}}{\sqrt{\gp e}}\norm{\gm}_{\mathfrak M(\prt Q_{\BBR_+}^T)}.
$$
This proves that 
\bel{IIl-1-7}
\norm{v_{0,\gm}\lfloor_{\prt Q_{\BBR_+}^T}}_{M^2(\prt Q_{\BBR_+}^T)}\leq c\norm{\gm}_{\mathfrak M(\prt Q_{\BBR_+}^T)}.
\ee
Similarly, if $G\subset [0,T]\ti[0,\infty)$ is a Borel set, then 
\bel{IIl-1-8}
\myint{G}{}\tilde E(t-s,x,0)ds\leq \myfrac{2|G|^{\frac 13}}{\sqrt{\gp e}},
\ee
and 
\bel{IIl-1-9}
\norm{v_{0,\gm}}_{M^3( Q_{\BBR_+}^T)}\leq c\norm{\gm}_{\mathfrak M(\prt Q_{\BBR_+}^T)}.
\ee
In the same way we prove that 
\bel{IIl-1-10}
\norm{v_{\gn,0}}_{M^3(Q^T_{\BBR_+})}+\norm{v_{\gn,0}\lfloor_{\prt Q_{\BBR_+}^T}}_{M^2(\prt Q_{\BBR_+}^T)}\leq c
\norm\gn_{\mathfrak M(Q^T_{\BBR_+})}.
\ee
This ends the proof.\qeda\medskip

\noindent{\it Proof of \rth{exist-uniq}}\smallskip

\noindent{\it Uniqueness.} Assume $u$ and $\tilde u$ are solutions of $(\ref{I-1-0})$, then $w=u-\tilde u$ satisfies 
\bel{III-1-1}\BA {lll}
\phantom{g(u---}
\phantom{,,,,-g(\tilde u(.,0))}w_t-w_{xx}=0\qquad&\text{in }\;Q^T_{\BBR_+}\\[1mm]
\phantom{;}-w_x(.,0)+g(u(.,0))-g(\tilde u(.,0))=0\qquad&\text{in }\;[0,T)\\[1mm]\phantom{-------,,-g(\tilde u(.,0))}
w(0,.)=0\qquad&\text{in }\;\BBR_+.
\EA
\ee
Applying $(\ref{II-1-9})$, we obtain
$$-\myint{0}{T}\!\!\myint{0}{\infty}(\gz_t+\gz_{xx})|w|dxdt+\myint{0}{\infty}(g(u(.,0))-g(\tilde u(.,0)) )sign (w)\gz(t,0) dt\leq 0,
$$
for any $\gz\in \BBX_{\BBR_+}^T$ with $\gz\geq 0$.  
Let $\gth\in C^1_c(Q^T_{\BBR_+})$, $\eta\geq 0$, we take $\gz$ to be the solution of
$$\BA {lll}-\gz_t-\gz_{xx}=\gth&\qquad\text{in }(0,T)\ti\BBR_+\\
\phantom{-\gz_,}
\gz_x(t,0)=0&\qquad\text{in }(0,T)\\
\phantom{-\gz_t}
\gz(T,x)=0&\qquad\text{in }(0,\infty).
\EA$$
Then $\gz$ satisfies $(\ref{II-1-14})$, hence
$$\myint{0}{T}\!\!\myint{0}{\infty}\gth|w|dxdt+\myint{0}{\infty}(g(u(.,0))-g(\tilde u(.,0)) )sign (w)\gz(t,0) dt\leq 0.
$$
This implies $w=0$. \smallskip

\noindent{\it Existence.} Without loss of generality we can assume that $\gm$ and $\gn$ are nonnegative. Let $\{\gn_n\}\subset C_c(\BBR_+)$ and $\{\gm_n\}\subset C_c([\BBR_+]0,T))$ converging to $\gn$ and $\gm$
 in the sense of measures and let $u_n$ be the solution of $(\ref{II-1-5})$. Then from $(\ref{II-1-7})$,
 \bel{III-1-2}\BA {lll}
 \myint{0}{T} \myint{0}{\infty}|u_n|dxdt+ \myint{0}{T}|g(u_n(t,0))|dt\leq T\myint{0}{\infty}|\gn_n|dx+ \myint{0}{T} |\gm_n|dt.
\EA
\ee
Therefore $u_n$ and $g(u_n(.,0))$ remain bounded respectively in $L^1(Q^T_{\BBR_+})$ and in $L^1(0,T)$. Furthermore, by \rlemma{reg1},
 $u_n$ remains bounded in $M^3(Q^T_{\BBR_+})$ and in $M^2(\prt Q^T_{\BBR_+})$. We can also write $u_n$ under the form
  \bel{III-1-3}\BA {lll}
  u_n(t,x)=\myint{0}{\infty}\tilde E(t,x,y)\gm_n(y)dy+2\myint{0}{t}E(t-s,x)(\gn_n(t)-g(u_n(t,0)))ds\\[2mm]
  \phantom{  u_n(t,x)}=A_n(t,x)+B_n(t,x).
 \EA
\ee
Since we can perform the even reflexion through $y=0$, the mapping 
$$(t,x)\mapsto A_n(t,x):=\myint{0}{\infty}\tilde E(t,x,y)\gm_n(y)dy,
$$
is relatively compact in $C^m_{loc}(\overline{Q^T_{\BBR_+}})$ for any $m\in\BBN^*$. Hence we can extract a subsequence $\{u_{n_k}\}$ which converges uniformly on every compact subset of $(0,T]\ti[0,\infty)$, hence a.e. on $(0,T]$ for the 1-dimensional Lebesque measure. Concerning the boundary term
$$(t,x)\mapsto B_n(t,x):=\myint{0}{t}E(t-s,x)(\gn_n(t)-g(u_n(t,0)))ds,
$$
it is relatively compact on every compact subset of $[0,T]\ti(0,\infty)$. If $x=0$, then 
$$B_n(t,0)=\myint{0}{t}(\gn_n(t)-g(u_n(t,0)))\myfrac{ds}{\sqrt{\gp(t-s)}}.
$$
Since $\norm{\gn_n(.)-g(u_n(.,0))}_{L^1(0,T)}$, $t\mapsto B_n(t,0)$ is uniformly integrable on $(0,T)$, hence relatively compact by the Frechet-Kolmogorov Theorem. Therefore there exists a subsequence, still denoted by $\{n_k\}$ such that $B_{n_k}(t,0)$ converges for almost all $t\in (0,T)$. This implies that the 
sequence of function $\{u_{n_k}\}$ defined by $(\ref{III-1-3})$ converges  in $\overline {Q^T_{\BBR_+}}$ up to a set $\Gth\cup\Gl$ where $\Gth\subset Q^T_{\BBR_+}$ is neglectable for the 2-dimensional Lebesgue measure and $\Gl\subset \prt_\ell Q^T_{\BBR_+}$ neglectable for the 1-dimensional Lebesgue measure.\par

From \rlemma{reg1}, $(u_{n,k}\lfloor_{Q^T_{\BBR_+}},u\lfloor_{\prt_\ell Q^T_{\BBR_+}})$ converges in $L^1_{loc}(Q^T_{\BBR_+})\ti L^1(\prt_\ell Q^T_{\BBR_+})$ and the convergence of each of the components holds also almost everywhere (up to a subsequence). Since $u_{n,k}$ is a weak solution, it satisfies for any $\gz\in\BBX(Q^T_{\BBR_+})$
 \bel{III-1-4}\BA {lll}
-\myint{0}{T}\!\!\myint{0}{\infty}(\gz_t+\gz_{xx})u_{n,k}dxdt+\myint{0}{T}\left(g(u_{n,k})\gz\right)(t,0)dt\\[4mm]
\phantom{---------}=\myint{0}{\infty}\gz \gn_{n,k}(x)dx+
\myint{0}{T}\gz (t,0) \gm_{n,k}(t)dt.
\EA
\ee
In order to prove the convergence of $g(u_{n,k}(t,0))$, we use Vitali's convergence theorem and the assumption $(\ref{I-1-15})$. Let $F\subset [0,T]$ be a Borel set. Using the fact that $0\leq u_{n,k}\leq v_{\gn_{n,k},\gm_{n,k}}$ and the estimate of \rlemma{reg1}, we have for any $\gl>0$, 
$$\BA {lll}
\myint{F}{}|g(u_{n,k}(t,0))|dt\leq \myint{F\cap\{u_{n,k}(t,0)\leq\gl\}}{}|g(u_{n,k}(t,0))|dt\\[4mm]
\phantom{---------\myint{F}{}|g(u_{n,k}(t,0))|dt}
+\myint{\{u_{n,k}(t,0)>\gl\}}{}|g(u_{n,k}(t,0))|dt\\[4mm]
\phantom{\myint{F}{}|g(u_{n,k}(t,0))|dt}
\leq g(\gl)|F|-\myint{\gl}{\infty}\gs d|\{t:|g(u_{n,k}(t,0))|>\gs\}|\\[4mm]
\phantom{\myint{E}{}|g(u_{n,k}(t,0))|dt}
\leq g(\gl)|F|+c\myint{\gl}{\infty}|g(\gs)|\gs^{-3}ds,
\EA$$
where $c$ depends of $\norm\gm_{\mathfrak M(\prt Q_{\BBR_+}^T)}+\norm\gn_{\mathfrak M(Q^T_{\BBR_+})}$. For $\ge>0$ given, we chose $\gl$ large enough so that the integral term above is smaller than $\ge$ and then $|F|$ such that $g(\gl)|F|+\leq\ge$. Hence $\{g(u_{n,k}(.,0))\}$ is uniformly integrable. 
Therefore up to a subsequence, it converges to $g(u(.,0))$ in $L^1(0,T)$. Clearly $u$ satisfies 
 \bel{III-1-5}\BA {lll}
-\myint{0}{T}\!\!\myint{0}{\infty}(\gz_t+\gz_{xx})udxdt+\myint{0}{T}\left(g(u)\gz\right)(t,0)dt\\[4mm]
\phantom{---------}=\myint{0}{\infty}\gz \gn(x)dx+
\myint{0}{T}\gz (t,0) \gm(t)dt,
\EA
\ee
which ends the existence proof.
\smallskip

\noindent{\it Monotonicity.} If $\gn\geq\tilde\gn$ and $\gm\geq\tilde\gm$; we can choose the approximations such that $\gn_n\geq\tilde\gn_n$ and $\gm_n\geq\tilde\gm_n$. It follows from $(\ref{II-1-10})$ that $u_{\gn_n,\gm_n}\geq u_{\tilde\gn_n,\tilde\gm_n}$. Choosing the same subsequence $\{n_k\}$, the limits
$u$, $\tilde u$ are in the same order. The conclusion follows by uniqueness. 
\qeda
\subsection{The case $g(u)=|u|^{p-1}u$}

Condition $(\ref{I-1-15})$ is satisfied if $p<2$. If this condition holds there exists a solution $u_{\ell\gd_0}=u_{0,\ell\gd_0}$ and the mapping 
$\ell\mapsto u_{\ell\gd_0}$ is increasing. 

\bth{lim-k} (i) If $1<p\leq \frac{3}{2}$, $u_{\ell\gd_0}$ tends to $\infty$ when $k\to\infty$.\smallskip

\noindent(ii) If $\frac{3}{2}<p<2$, $u_{\ell\gd_0}$ converges to $U_{\gw_s}$ defined by 
$$U_{\gw_s}(t,x)=t^{-\frac{1}{2(p-1)}}\gw_s(\tfrac{x}{\sqrt t}),$$ 
when $k\to\infty$.
\es
\Proof By uniqueness and using $(\ref{I-1-3})$, there holds
 \bel{III-3-1}\BA {lll}
T_k[u_{\ell\gd_0}]=u_{k^{\frac{2-p}{p-1}\ell }\gd_0},
\EA
\ee
for any $k,\ell>0$. Since $\ell\mapsto u_{\ell\gd_0}$ is increasing, its limit $u_\infty$, when $\ell\to\infty$, satisfies 
 \bel{III-3-2}\BA {lll}
T_k[u_{\infty}]=u_{\infty}. 
\EA
\ee
Hence $u_\infty$ is a positive self-similar solution of $(\ref{I-1-2})$, provided it exists. Hence $u_\infty=U_{\gw_s}$ if $\frac{3}{2}<p<2$. If 
 $1<p\leq \frac{3}{2}$, $u_{k\gd_0}$ admits no finite limit when $k\to\infty$ which ends the proof. \qeda \medskip
 
 \noindent\Remark As a consequence of this result, no a priori estimate of Brezis-Friedman type (parabolic Keller-Osserman) exists for a nonnegative function 
 $u\in C^{2,1}(\overline{Q^\infty_{\BBR_+}}\setminus\{(0,0)\}$
solution of
 \bel{III-3-3}\BA {lll}
\phantom{-------,,ii}
u_t-u_{xx}=0\quad&\text{in }\;Q^\infty_{\BBR_+}\\[0mm]\phantom{--}
-u_x(.,0)+|u|^{p-1}u(.,0)=0&\text{for all }\;t>0
\\[1mm]
\phantom{,,--------,,i}\displaystyle
u(0,x)=0&\text{for all }\;x>0.
\EA
\ee
when $1<p\leq \frac{3}{2}$. When $\frac{3}{2}<p<2$ it is expected that 
 \bel{III-3-4}\BA {lll}
u(t,x)\leq \myfrac{c}{(|x|^2+t)^{\frac{1}{2(p-1)}}}.
\EA
\ee
The type of phenomenon (i) in \rth{lim-k} is characteristic of fractional diffusion. It has already been observed in \cite[Theorem 1.3]{CVW} with equations
 \bel{III-3-5}\BA {lll}
u_t+(-\Gd)^{\ga}u+t^{\gb}u^p=0\quad &\text{ in }\, \BBR_+\ti\BBR^N\\
\phantom{-----\,\;\;-}
u(0,.)=k\gd_0&\text{ in }\, \BBR^N,
\EA
\ee
when $0<\ga<1$ is small and $p>1$ is close to $1$.

\section{Extension and open problems}
The natural extension is to replace a one dimensional domain by a mutidimenional one. The main open problem is the question of a priori estimate as stated in the last remark above. 
\subsection{Self-similar solutions}
Let $\eta=(\eta_1,...,\eta_n)$ be the coordinates in $\BBR^n$ and denote $\BBR_+^n=\{\eta=(\eta_1,...,\eta_n)=(\eta',\eta_n):\eta_n>0\}$. 
We set $K(\eta)=e^{\frac{|\eta|^2}{4}}$ and $K'(\eta')=e^{\frac{|\eta'|^2}{4}}$. 
Similarly to Section 2  we define $\CL_K$ in $C_0^2(\BBR^n)$ by 
\bel{IV-1-1}
\BA{lll}
\CL_K(\phi)=-K^{-1}div (K\nabla\phi).
\EA
\ee
If $\ga=(\ga_1,...,\ga_n)\in\BBN^n$, we set $|\ga|=\ga_1+\ga_2+...+\ga_n$. We denote by $\phi_1$ the function $K^{-1}$. Then the set of eigenvalues of $\CL_K$ is the set of numbers $\left\{\gl_k=\frac{n+k}{2}:k\in\BBN\right\}$ with corresponding set of eigenspaces
$$N_k=span\left\{ D^\ga\phi_1:|\ga|=k\right\}.
$$
The operators $\CL_K^{+,N}$ and $\CL_K^{+,D}$ are defined acoordingly in $H_K^{1}(\BBR^n_+)$ and $H_K^{1,0}(\BBR^n_+)$ respectively and $\gs(\CL_K^{+,N})=\left\{\frac{n+k}{2}:k\in\BBN\right\}$ and $\gs(\CL_K^{+,D})=\left\{\frac{n+k}{2}:k\in\BBN^*\right\}$
Furthermore
\bel{IV-1-2}N_{k,N}=ker\left(\CL_K^{+,N}-\tfrac{n+k}{2}I_d\right)=span\left\{ D^\ga\phi_1:|\ga|=k,\ga_n=2\ell\,,\;\ell\in\BBN\right\},
\ee
and
\bel{IV-1-3}N_{k,D}=ker\left(\CL_K^{+,D}-\tfrac{n+k}{2}I_d\right)=span\left\{ D^\ga\phi_1:|\ga|=k,\ga_n=2\ell+1\,,\;\ell\in\BBN\right\}.
\ee
Since $\CL_K^{+,N}$ and $\CL_K^{+,D}$ are Fredholm operators, 
\bel{IV-1-4}
H^1_K(\BBR^n_+)=\bigoplus_{k=0}^\infty N_{k,N}\;\text{ and }\;H^{1,0}_K(\BBR^n_+)=\bigoplus_{k=1}^\infty N_{k,D}.
\ee
 We define the following functional on $H^1_K(\BBR^n_+)$
   \bel{IV-1-5}
   J(\phi)=\myfrac{1}{2}\myint{\BBR^n_+}{}\left(|\nabla\gf|^2-\myfrac{1}{2(p-1)}\gf^2\right) Kd\eta+\myfrac{1}{p+1}\myint{\prt\BBR^n_+}{}|\gf|^{p+1}
   K'd\eta'.
\ee
The critical points of $J$ satisfies 
 \bel{IV-1-6}\BA {lll}
\!\!-\Gd\gw-\myfrac12\eta.\nabla\gw-\myfrac1{2(p-1)}\gw=0\quad&\text{in }\; \BBR^n_+\\[2mm]
\phantom{---,-\prt}
-\gw _{\eta_n}+|\gw|^{p-1}\gw=0\quad&\text{in }\; \prt\BBR^n_+.\\[1mm]
\EA
\ee 
If $\gw$ is a solution of $(\ref{IV-1-6})$, the function
 \bel{IV-1-7}
 u_\gw(t,x)=t^{-\frac{1}{2(p-1)}}\gw(\frac x{\sqrt t})
 \ee
 satisfies
  \bel{IV-1-8}\BA {lll}\phantom{----}
u_{\gw\,t}-\Gd u_\gw=0\qquad&\text{in }Q_{\BBR^n_+}^\infty:=(0,\infty)\ti\BBR^n_+\\[1mm]
-u_{\gw\,x_n}+|u_\gw|^{p-1}u_\gw=0\qquad&\text{in }\prt_\ell Q_{\BBR^n_+}^\infty:=(0,\infty)\ti\prt\BBR^n_+.
 \EA\ee
 Here we have set $\BBR_+^n=\{x=(x_1,...,x_n)=(x',x_n):x_n>0\}$. We denote by $\CE$ the subset  $H^{1}_K(\BBR^n_+)\cap L^p(\prt\BBR^n_+;d\eta')$ of solutions of $(\ref{IV-1-6})$ and by $\CE_+$ the subset of positive solutions.
 As for the case $n=1$ we have the following non-existence result
 \bprop{nonex} 1- If $p\geq1+\frac{1}{n}$, then $\CE=\{0\}$.\smallskip

\noindent 2- If $1<p\leq 1+\frac{1}{n+1}$, then $\CE_+=\{0\}$\smallskip
 \es
 The proof is similar to the one of \rth{vssth}. Hence the existence is to be found in the range $1+\frac{1}{n+1}<p<1+\frac 1n$. \medskip

 \noindent {\bf Conjecture} {\it Assume $1+\frac{1}{n+1}<p<1+\frac 1n$, then the functional $J$ is bounded from below in $H^{1}_K(\BBR^n_+)\cap L_{K'}^p(\prt\BBR^n_+)$. Furthermore $J(\phi)$ tends to infinity when $\norm {\phi}_{H^{1}_K(\BBR^n_+)}+\norm {\phi\lfloor_{\prt\BBR^n_+}}_{L^{p+1}_{K'}(\prt\BBR^n_+)}$ tends to infinity.}

 \subsection{Problem with measure data}
 The method for proving \rth{exist-uniq} can be adapted to prove the following $n$-dimensional result
 \bth{exist-uniq-n} Let $g:\BBR\mapsto\BBR$ be a nondecreasing continuous function such that $g(0)=0$ and 
   \bel{IV-1-10}\BA {lll}\phantom{----}
\myint{1}{\infty}(g(s)-g(-s))s^{-\frac{2n+1}{n}}ds<\infty,
 \EA\ee
 then for any bounded Radon measures $\gn$ in $\BBR^n_+$ and $\gm$ in $(0,T)\ti\prt\BBR^n_+$, there exists a unique Borel function $u:=u_{\gn,\gm}$ defined 
 in $\overline{Q_{T}^{\BBR^n_+}}:=[0,T]\ti\BBR^n_+$ such that $u\in L^1(Q_{T}^{\BBR^n_+})$, $u\lfloor_{(0,T)\ti\prt\BBR^n_+}\in L^1((0,T)\ti\prt\BBR^n_+)$ and $g(u)\in L^1((0,T)\ti\prt\BBR^n_+)$ solution of
 \bel{IV-1-11}\BA {lll}
\phantom{g(u;}
\phantom{,--,}u_t-\Gd u=0\qquad&\text{in }\;Q^T_{\BBR^n_+}\\[1mm]\phantom{,,,u}
-u_{x_n}+g(u)=\gm\qquad&\text{in }\;\prt_\ell Q^T_{\BBR^n_+}\\[1mm]\phantom{--,,---}
u(0,.)=\gn\qquad&\text{in }\;\BBR^n_+,
\EA
\ee
in the sense that 
 \bel{IV-1-12}\BA {lll}
\myint{}{}\myint{Q^T_{\BBR^n_+}}{}(-\prt_t\gz-\Gd\gz)udxdt+\myint{}{}\myint{\prt_\ell Q^T_{\BBR^n_+}}{}g(u)\gz dx'dt\\[4mm]
\phantom{----------------}=\myint{\BBR^n_+}{}\gz d\gn+
\myint{}{}\myint{\prt_\ell Q^T_{\BBR^n_+}}{}\gz d\gm,
\EA
\ee
for all $\gz\in C_c^{1,2}(\overline{Q^T_{\BBR^n_+}})$ such that $\gz_{x_n}=0$ on $(0,T)\ti\prt\BBR^n_+$ and $\gz(T,.)=0$. Furthermore $(\gn,\gm)\mapsto u_{\gn,\gm})$ is nondecreasing. 
 \es
\bibliographystyle{Bibtex}

\begin{thebibliography}{110}

\bibitem{AS} M. A{\scriptsize{BRAMOWITZ}}, I. A. S{\scriptsize{TEGUN}}. \textit{Handbook of Mathematical Functions}, National Bureau of Standards,
Washington, 1964.

\bibitem{BoV} O. B{\scriptsize{OUKARABILA}}, L. V{\scriptsize{\'ERON}}. \textit{Nonlinear boundary value problems \\relative to harmonic functions}, Nonlinear Analysis, to appear. arXiv:2003.00871. 

\bibitem{Br-OMM} H. B{\scriptsize{REZIS}}. \textit{ Op\'erateurs maximaux monotones et semi-groupes de contractions dans les espaces de Hilbert}, Notas de Matem\`aticas {\bf 5}, North Holland (1971).

\bibitem{BrFr} H. B{\scriptsize{REZIS}}, A. F{\scriptsize{RIEDMAN}}. \textit{Nonlinear parabolic equations involving measures as initial
conditions}, J. Math. Pures Appl. {\bf 62} (1983), 73--97.

\bibitem{BPT} H. B{\scriptsize{REZIS}}, L.A. P{\scriptsize{ELETIER}}, D. T{\scriptsize{ERMAN}}. \textit{A very singular solution of the heat equation with
absorption}, Arch. Rational Mech. Anal. {\bf 95} (1986), 185--209.

\bibitem{CVW} H. C{\scriptsize{HEN}}, L. V{\scriptsize{\'ERON}}, Y. W{\scriptsize{ANG}}. \textit{Fractional heat equations with subcritical absorption having a
measure as initial data}, Nonlinear Anal. {\bf 137} (2016), 306--337.

\bibitem{EK} M. E{\scriptsize{SCOBEDO}}, O. K{\scriptsize{AVIAN}}. \textit{Variational problems related to the self-similar solutions of the heat
equations}, Nonlinear Anal. {\bf 10} (1987), 1103--1133. 

\bibitem {Fil} M. F{\scriptsize{ILA}}, K. I{\scriptsize{SHIGE}}, T. K{\scriptsize{AWAKAMI}}. \textit{Existence of positive solutions of a semilinear elliptic
equation with a dynamical boundary condition}, Calc. Var. {\bf 54} (2015), 2059--2078.


\bibitem{FQ} M. F{\scriptsize{ILA}}, P. Q{\scriptsize{UITTNER}}. \textit{The blow-up rate for the heat equation with a non-linear
boundary condition}, Math. Methods in the Appl. Sci. {\bf 14} (1991), 197--205.

\bibitem {Ish-Ka} K. I{\scriptsize{SHIGE}}, T. K{\scriptsize{AWAKAMI}}. \textit{Global solutions of the heat equation with a nonlinear boundary condition}, Calc. Var. {\bf 39} (2010), 429--457.

\bibitem {Ish-Sa} K. I{\scriptsize{SHIGE}}, R. S{\scriptsize{ATO}}. \textit{Heat equation with a nonlinear boundary condition and
uniformly local $L^r$ spaces}, Disc. Cont. Dyn. Syst. {\bf 36} (2016), 2627--2652.

\bibitem{Ke} J. K{\scriptsize{ELLER}}. \textit{On solutions of $\Delta u = f(u)$}, Comm. Pure Appl. Math. \textbf{10} (1957), 503--510.

\bibitem{MV1} M. M{\scriptsize{ARCUS}}, L. V{\scriptsize{\'ERON}}. \textit{Semilinear parabolic equations with measure boundary data and isolated singularities}, J. Anal. Mat.  {\bf 85}  (2001), 245--290.

\bibitem{MV2} M. M{\scriptsize{ARCUS}}, L. V{\scriptsize{ERON}}. \textit{Isolated boundary singularities of signed solutions of some nonlinear parabolic equations}, Adv. Diff. Equ. {\bf 6} (2001), 1281--1316.

\bibitem{MoV} I. M{\scriptsize{OUTOUSSAMY}}, L. V{\scriptsize{\'ERON}}. \textit{Isolated singularities and asymptotic behaviour of the solutions of a semilinear heat equation}, Asymptotic Anal. {\bf 9}  (1994), 259--289.

\bibitem{Os} R. O{\scriptsize{SSERMAN}}. \textit{On the inequality $\Delta u \geq f(u)$}, Pacific J. Math. \textbf{7} (1957), 1641--1647.
\end{thebibliography}



\authinfo{ 
Laurent V\'eron\\
Institut Denis Poisson, CNRS UMR 7013, \\
 Universit\'e de Tours, France.\\
\email{veronl@univ-tours.fr}}

\end{document}